\documentclass[12pt]{article}
\usepackage[mathscr]{eucal}
\usepackage{amssymb,amsmath,amscd}
\font\cmc=cmcsc10  scaled \magstep2

\newcommand\Q{\mathbb{Q}}

\newcommand\mscrK{\mathscr{K}}

\newcommand\vk{\vskip}
\newcommand\hk{\hskip}
\newcommand\al{\alpha}
\newcommand\be{\beta}
\newcommand\iy{\infty}

\newcommand\ga{\gamma}

\newcommand\sg{\sigma}

\newcommand\rg{\rightarrow}

\newcommand\om{\omega}

\newcommand\var{\varphi}

\newcommand\no{\noindent}

\newcommand\cen{\centerline}

\newtheorem{Proof.}{\it Proof.}

\begin{document}
\vbox to .5truecm{}

\begin{center}
\cmc 
A Genus Two Curve Related\\
to the Class Number One Problem
\end{center}
\vk.3cm
\begin{center}
by 
Viet Kh. Nguyen
\footnote{IMC Institute \& Fpt University\\
{\text {\hk.5cm}} 2010 Mathematics Subject Classification. Primary 11G30, 11R29, 11R11.}
\end{center}
\vk.1cm
\hk10.5cm{\it \`a la m\'emoire de V. A. Iskovskikh}

\vk.2cm

\begin{abstract} {We give another solution to the class number one problem by showing that imaginary quadratic fields $\Q(\sqrt{-d})$ with class number $h(-d)=1$ correspond to integral points on a genus two curve $\mscrK_3$. In fact one can find all rational points on $\mscrK_3$. The curve $\mscrK_3$ arises naturally via certain coverings of curves:\ $\mscrK_3\rg\mscrK_6$,\ $\mscrK_1\rg\mscrK_2$\ with $\mscrK_2\colon y^2=2x(x^3-1)$ denoting the Heegner curve, also in connection with the so-called  Heegner-Stark covering $\mscrK_1\rg\mscrK_s$.}
\end{abstract}

\vk.5cm					
{\bf 1.\ Introduction.}\quad In his famous {\it Disquisitiones Arithmeticae} (Article 303) Gauss hypothetically believed that there are only nine 
imaginary quadratic fields $\Q(\sqrt{-d})$ with class number $h(-d)=1$,\ namely with\ $d=1, 2, 3, 7, 11, 19, 43, 67, 163$. It is known nowadays as 
the Baker-Heegner-Stark theorem. 
\vk.2cm
In the early 1980's V. A. Iskovskikh organized a seminar devoted to a beautiful approach of Heegner who was substantially first to give a 
complete proof of the theorem of Baker-Heegner-Stark based on an approach using modular functions (\cite{Heeg52}, {\it cf.} 
also \cite{Bir69}, \cite{Star73}, \cite{Ser89}, \cite{Star07}, \cite{Sche10}, \cite{Sha13}). More precisely, assuming $d\equiv 3\ ({\rm mod}\ 8)$ 
(the main case), Heegner showed that every such a field corresponds to an integral point on the so-called Heegner curve:   
  $$y^2=2x(x^3-1)\leqno{(\mscrK_2)}$$

In the higher class number case the statement reads as every field $\Q(\sqrt{-d})$,\ $d\equiv 3\ ({\rm mod}\ 8)$ corresponds to an integral point of 
degree $h(-d)$ on $\mscrK_2$, provided $3\nmid d$. One may have a remark that the ``{\it exceptional}" field $\Q(\sqrt{-3})$ (corresponding to the 
point $(0,0)$ on $\mscrK_2$) is the only case with $3\mid d$ such that $j$-invariant is a cube (namely equal to zero). On the contrary to the case of 
higher class numbers one sees, for example if $h(-d)=2$, three fields 
$\Q(\sqrt{-51})$,\ $\Q(\sqrt{-123})$ and $\Q(\sqrt{-267})$ have no representatives on $\mscrK_2$ along this line of arguments.    
\vk.2cm
The aim of this note is to show that for our purpose one can involve another genus two curve $\mscrK_3$:
  $$8 x^8-32 x^6 y+40 x^4 y^2+64 x^5-16 x^2 y^3-128 x^3 y+y^4+48 x y^2+96 x^2-32 y-24=0\leqno{(\mscrK_3)}$$
realized at that time by the author. 

The main result of the note is  
\vk.2cm
{\bf Theorem A.}\quad {\it The curve $\mscrK_3$ has the following integral points 
\vk.2cm
\cen {\quad\ $(3,6)$,\quad\ $(-1,2)$,\qquad  $(1,6)$,\qquad $(3,14)$,\qquad $(7,26)$,\quad $(-17,150)$ }
\vk.2cm
\no corresponding to the fields
\vk.2cm
\cen{ $\Q(\sqrt{-3})$,\ $\Q(\sqrt{-11})$,\ $\Q(\sqrt{-19})$,\ $\Q(\sqrt{-43})$,\ $\Q(\sqrt{-67})$,\ 
$\Q(\sqrt{-163})$}
\vk.2cm
\no respectively. Besides it has three ``extraneous" integral points:\ $(-1,-2)$,\ $(-3,6)$,\ $(1,2)$\ (in fact $(1,2)$\ is a double 
point on $\mscrK_3$)\ plus two more rational points\ $\big(-\frac{9}{17},\frac{6}{289}\big)$,\ $\big(-\frac{155}{79},\frac{42486}{6241}\big)$.   
All the listed are the only rational points on $\mscrK_3$.}
\vk.2cm
As a corollary one obtains a new solution to Gauss' ``{\it tenth discriminant problem}". 
\vk.2cm
Nowadays with the aid of computational tools one can easily compute the genus of $\mscrK_3$. 
A surprising observation is the fact that $\mscrK_3$ is birational to a genus 2 curve $\mscrK_s$ studied independently by Heegner and Stark 
(\cite{Heeg52}, \cite{Star73}). 
\vk.5cm
{\bf Acknowledgments.}\quad I would like to thank V. A. Abrashkin for a nice talk on the class number problem at the seminar mentioned in the 
Introduction. This paper is dedicated to the 75th birthday of my teacher V. A. Iskovskikh with profound admiration and memory.  
\vk.5cm
{\bf 2.\ The curve $\mscrK_3$.}\quad Let us first recall some known results from the Weber-Heegner-Birch theory (\cite{Web91}, 
\cite{Heeg52}, \cite{Bir69}, \cite{Sche10}). We shall use Heegner's notation $\sg$ for the Schl\"afli functions of Stufe $48$ which are defined by 
  $$\sg(\tau)=q^{-\frac{1}{24}}\ \prod_{n=1}^\iy\ \big( 1+q^{2n-1}\big);\quad 
\sg_1(\tau)=q^{-\frac{1}{24}}\ \prod_{n=1}^\iy\ \big( 1-q^{2n-1}\big);$$
  $$\sg_2(\tau)=\sqrt{2}\ q^{\frac{1}{12}}\ \prod_{n=1}^\iy\ \big( 1+q^{2n}\big);$$
with $q=e^{i\pi\tau}$,\ Im$(\tau)>0$. Then the $24th$ powers $\sg^{24}$,\ $-\sg_1^{24}$,\ $-\sg_2^{24}$ are the roots of
  $$U^3-48 U^2+(768-j)U-4096=0\eqno{(2.1)}$$
where $j$ is the modular invariant.
\vk.2cm
We assume throughout $d$ is positive $\equiv 3\ ({\rm mod}\ 8)$ and put $\om=\frac{1}{2}(3+\sqrt{-d})$. Then $-\sg_2^{24}(\om)$ is real and positive. 
Let $\phi$ denote the $3rd$ power $\sg_2^3(\om)$ up to a suitable $16th$ root of unity such that $\phi$ is real and positive. Then $\phi^8$ is an 
algebraic integer satisfying equation $(2.1)$ with $j=j(\om)$ and by Weber's results $\phi^4,\ \phi^2$ are algebraic integers of degree $3h(-d)$. 
Furthermore Weber's conjecture proved in \cite{Bir69} tells us that indeed $\sqrt{2}\ \phi$ is integral of degree $3h(-d)$, so it satisfies a cubic 
equation
  $$W^3-2 a_3 W^2+2 b_3 W-8=0\eqno{(2.2)}$$
with $a_3, b_3$ integral of degree $h(-d)$. A computation in conjunction with $(2.1)$ shows that $(a_3,b_3)$ gives rise to an integral point of 
degree $h(-d)$ on $\mscrK_3$. In particular if $h(-d)=1$ we get the first statement of Theorem A.
\vk.2cm
If we argue in the same way with a cubic equation for $\phi^2$
  $$T^3-2 a_2 T^2+2 b_2 T-8=0\eqno{(2.3)}$$
we get a relation between $a_2, b_2$ describing a curve $\mscrK_6$ which is birationally isomorphic to the Pell conics
  $$\Big(\frac{k}{2}-(a_2+1)\Big)^2-2\Big(\frac{a_2+1}{2}\Big)^2=1$$
where\ $b_2=k(a_2-1)+2$. 
\vk.2cm
In view of $(2.2)$,\ $(2.3)$ the covering $\mscrK_3\rg\mscrK_6$ can be given in the following explicit form
  $$a_2=a_3^2-b_3,\quad 2 b_2=b_3^2-8 a_3\eqno{(2.4)}$$
from which may find interesting relations with the Pell sequence.
\vk.2cm
{\bf Remark\ 2.1.}\quad By virtue of $(2.4)$ one sees that $a_3$ satisfies the following equation
  $$z^4-2 a_2 z^2-8z+(a_2^2-2 b_2)=0\eqno{(2.5)}$$
Therefore $(2.3)$ is nothing but Euler's resolvent for $(2.5)$.
\vk.2cm
{\bf Example\ 2.2.}\quad As noted in the Introduction we illustrate our approach with some values $d$ multiple of $3$. By using the Tables from 
\cite{Web91} and \cite{Wat36} one sees that the fields $\Q(\sqrt{-51})$,\ $\Q(\sqrt{-123})$,\ $\Q(\sqrt{-267})$ give rise to three real quadratic 
points 
  $$(-1,8+2\sqrt{17}),\quad (4+\sqrt{41},40+6\sqrt{41}),\quad (-10-\sqrt{89},310+32\sqrt{89})$$
on $\mscrK_3$ respectively. 
\vk.5cm
{\bf 3.\ The curve $\mscrK_s$.}\quad If we assume in addition $3\nmid d$, then the $8th$ powers 
$\sg^8$,\ $-\sg_1^8$,\ $-\sg_2^8$ are the roots of
  $$V^3-\ga_2 V-16=0\eqno{(3.1)}$$
where $\ga_2$ is the Weber function of Stufe $3$ such that $\ga_2^3=j$.
\vk.2cm 
In the notation of {\bf 2.} let $\var$ be a cubic root of $\phi$ that is real and positive. Then as above $\var^2$,\ $\sqrt{2}\ \var$ are algebraic 
integral of degree $3h(-d)$. From $(3.1)$ and a cubic equation for $\var^2$
  $$Z^3-2 \al_2 Z^2+2 \be_2 Z-2=0\eqno{(3.2)}$$
with $\al_2, \be_2$ integral of degree $h(-d)$ we come to the equation ($\mscrK_2$) of Heegner's curve with $x=\al_2$,\ $y=\be_2-2\al_2^2$. 
\vk.2cm
Now writing a cubic equation for $\sqrt{2}\ \var$
  $$S^3-2 \al_3 S^2+2 \be_3 T-4=0\eqno{(3.3)}$$
with $\al_3, \be_3$ integral of degree $h(-d)$ one gets a covering\ $\mscrK_1\rg\mscrK_2$ considered in \cite{Heeg52}, \cite{Star73} which can be 
given by 
  $$\al_2=\al_3^2-\be_3,\quad 2 \be_2=\be_3^2-4 \al_3\eqno{(3.4)}$$ 

Hence the defining equation for $\mscrK_1$ is as follows 
  $$8 x^8-32 x^6 y+40 x^4 y^2+32 x^5-16 x^2 y^3-64 x^3 y+y^4+24 x y^2+24 x^2-8 y=0\leqno{(\mscrK_1)}$$

Despite a similarity with equation $(\mscrK_3)$ this is a curve of genus $9$.
\vk.2cm
{\bf Proposition\ 3.1.}\quad {\it The curve $\mscrK_1$ has exactly 6 integral points $(0,0)$,\ $(1,2)$,\ $(-1,0)$,\ $(0,2)$,\ $(-1,2)$,\ $(2,6)$ 
which are in $1-1$ correspondence with class number one fields 
$\Q(\sqrt{-3})$,\ $\Q(\sqrt{-11})$,\ $\Q(\sqrt{-19})$,\ $\Q(\sqrt{-43})$,\ $\Q(\sqrt{-67})$,\ $\Q(\sqrt{-163})$.}
\vk.2cm
The proof follows immediately from $(3.4)$ and \cite{Heeg52}. The fact that these integral points exhaust all rational points on $\mscrK_1$ follows 
from Heegner-Stark results (\cite{Heeg52}, \cite{Star73}) on rational points on curve $\mscrK_s$
  $$w^2=2z(z^4+4 z^3-2 z^2+4z+1)\leqno{(\mscrK_s)}$$

The covering $\mscrK_1\rg\mscrK_s$ is given by
  $$z=\frac{\be_3}{\al_3^2}-1;\quad w=4(z-2)y-2(3 z^2-2 z-1),\quad y=\frac{1}{\al_3^3}$$

It should be noted that Heegner used a little different normalization ({\it cf.} \cite{Heeg52}, formulae $(151)-(152)$). Furthermore in order to 
enumerate all rational points of $\mscrK_s$ one may apply the argument of \cite{Heeg52} in solving equation $(\mscrK_2)$ in integers. A precise 
statement is as follows.
\vk.2cm
{\bf Theorem\ 3.2.}\quad {\it The only rational points on the curve $\mscrK_s$ are\
  $(0,0)$,\ $(\pm 4,1)$,\ $(\pm 4,-1)$,\ $\big(\pm\frac{7}{4},\frac{1}{2}\big)$,\ $(\pm 14,2)$.
 }
\vk.2cm
{\bf Theorem\ 3.3.}\quad {\it The curves $\mscrK_3$ and $\mscrK_s$ are birationally isomorphic.}
\vk.2cm
Since $\mscrK_3$ is hyperelliptic of genus $2$ one may use van Hoeij's algorithm (\cite{vHoe02}) and a computing program ({\it e.g.} Maple) to write 
down its Weierstrass normal form which turns out to be the equation $(\mscrK_s)$. The corresponding birational isomorphism $\mscrK_s\rg\mscrK_3$ is 
given by    
  $$x=-\frac{z^4+8 z^3+2 w z+18 z^2+6 w-3}{z^4+4 z^3-2 z^2-12 z+1}$$
  $$y=\frac{2\ P_8(z,w)}{(z^4+4 z^3-2 z^2-12 z+1)^2}$$
where 
  $$P_8(z,w)=(2 z^5+10 z^4+36 z^3+68 z^2+10 z-30)w$$
  $$\hk4.5cm -z^8+60 z^6+192 z^5+82 z^4-128 z^3+172 z^2+64 z+7.$$
\vk.2cm
Theorem A now follows from Theorems 3.2 and 3.3.
\vk.2cm 
We note that the inverse birational isomorphism $\mscrK_3\rg\mscrK_s$ has more complicated formulae for $(z,w)$. For completeness we include them 
here in the following form 
  $$z=1-\frac{4 x^3-4 x y-y^2+4 x+4}{2 x^4+2 x^3-3 x^2 y-2 x y+6 x-y+2}$$
  $$w=-\frac{2\ P_{12}(x,y)}{(x-1)(x^2+1)(x^2-2 x-1)(x^2+2 x+3)(x+1)^5}$$
where 
  $$\begin{aligned}
P_{12}(x,y)\ =\ \ & 4 x^{12}+252 x^{11}-24 x^{10} y+156 x^{10}-622 x^9 y+15 x^8 y^2+440 x^9-514 x^8 y+\\
             + & 322 x^7 y^2-x^6 y^3+1256 x^8-708 x^7 y+288 x^6 y^2-21 x^5 y^3+1536 x^7-620 x^6 y\\
             + & 310 x^5 y^2-19 x^4 y^3+1344 x^6-716 x^5 y+64 x^4 y^2-20 x^3 y^3+440 x^5-640 x^4 y\\
             + & 22 x^3 y^2-7 x^2 y^3-12 x^4-316 x^3 y-8 x^2 y^2-3 x y^3-124 x^3-140 x^2 y-6 x y^2\\
             - & y^3-92 x^2-22 x y+y^2-16 x+2 y.
\end{aligned}$$ 
\vk.2cm
{\bf Remark\ 3.4.}\quad There is a map $\mscrK_2\rg\mscrK_6$ (\cite{Heeg52}) given by
  $$a_2=4\al_2^3-6\al_2\be_2+3;\quad b_2=4\be_2^3-12\al_2\be_2+6$$
and similarly a map $\mscrK_1\rg\mscrK_3$ 
  $$a_3=2\al_3^3-3\al_3\be_3+3;\quad b_3=\be_3^3-6\al_3\be_3+6$$
which turned out quite useful in the computation process.
\vk.5cm

\vk.5cm
\no IMC Institute \& Fpt University 

\no Permanent: 8 Ton That Thuyet, My Dinh, Tu Liem, Hanoi, Vietnam

\no e-mail:\ vietnk@fpt.edu.vn,\ \ nhimsocvn@gmail.com

\end{document}